\documentclass[12pt]{article}%
\usepackage{graphicx}
\usepackage[intlimits]{amsmath}
\usepackage{amsfonts}
\usepackage{amssymb}%
 \makeatletter 
\newtheorem{conjecture}{Conjecture}

\newtheorem{theorem}{Theorem}

\begin{document}

\title{Cliques and the Spectral Radius}
\author{B\'{e}la Bollob\'{a}s\thanks{Department of Mathematical Sciences, University
of Memphis, Memphis TN 38152, USA} \thanks{Trinity College, Cambridge CB2 1TQ,
UK} \thanks{Research supported in part by DARPA grant F33615-01-C-1900.} \ and
Vladimir Nikiforov$^{\ast}$}
\maketitle

\begin{abstract}
We prove a number of relations between the number of cliques of a graph $G$
and the largest eigenvalue $\mu\left(  G\right)  $ of its adjacency matrix. In
particular, writing $k_{s}\left(  G\right)  $ for the number of $s$-cliques of
$G$, we show that, for all $r\geq2,$%
\[
\mu^{r+1}\left(  G\right)  \leq\left(  r+1\right)  k_{r+1}\left(  G\right)
+\sum_{s=2}^{r}\left(  s-1\right)  k_{s}\left(  G\right)  \mu^{r+1-s}\left(
G\right)  ,
\]
and, if $G$ is of order $n,$ then
\[
k_{r+1}\left(  G\right)  \geq\left(  \frac{\mu\left(  G\right)  }{n}%
-1+\frac{1}{r}\right)  \frac{r\left(  r-1\right)  }{r+1}\left(  \frac{n}%
{r}\right)  ^{r+1}.
\]

\textbf{Keywords: }\textit{number of cliques, clique number, spectral radius,
stability }

\end{abstract}

\section{Introduction}

Our graph-theoretic notation is standard (e.g., see \cite{Bol98}); in
particular, we write $G\left(  n\right)  $ for a graph of order $n$. Given a
graph $G,$ a $k$\emph{-walk} is a sequence of vertices $v_{1},\ldots,v_{k}$ of
$G$ such that $v_{i-1}$ is adjacent to $v_{i}$ for all $i=2,\ldots,k.$ We
write $w_{k}\left(  G\right)  $ for the number of $k$-walks in $G$ and
$k_{r}\left(  G\right)  $ for the number of its $r$-cliques. We order the
eigenvalues of the adjacency matrix of a graph $G=G\left(  n\right)  $ as
$\mu\left(  G\right)  =\mu_{1}\left(  G\right)  \geq\ldots\geq\mu_{n}\left(
G\right)  $.

Let $\omega=\omega\left(  G\right)  $ be the clique number of $G$. Wilf
\cite{Wil86} proved that
\[
\mu\left(  G\right)  \leq\frac{\omega-1}{\omega}v\left(  G\right)
=\frac{\omega-1}{\omega}w_{1}\left(  G\right)  ,
\]
and Nikiforov~\cite{Nik06} extended this, showing that the inequality%
\begin{equation}
\mu^{s}\left(  G\right)  \leq\frac{\omega-1}{\omega}w_{s}\left(  G\right)
\label{maxmu}%
\end{equation}
holds for every $s\geq1.$ Note that for $s=2$ inequality \eqref{maxmu} implies
a concise form of Tur\'{a}n's theorem. Indeed, if $G$ has $n$ vertices and $m$
edges, then $\mu(G)\geq2m/n,$ and so,
\[
\left(  \frac{2m}{n}\right)  ^{2}\leq\mu^{2}(G)\leq\frac{\omega-1}{\omega
}w_{2}(G)=\frac{\omega-1}{\omega}2m.
\]
This shows that
\begin{equation}
m\leq\frac{\omega-1}{2\omega}n^{2}, \label{maxmu1}%
\end{equation}
which is best possible whenever $\omega$ divides $n.$ If we combine
\eqref{maxmu} with other lower bounds on $\mu(G)$, e.g., with
\[
\mu^{2}(G)\geq\frac{1}{n}\sum_{u\in V(G)}d^{2}\left(  u\right)  ,
\]
we obtain generalizations of (\ref{maxmu1}).

Moreover, inequality \eqref{maxmu} follows from a result of Motzkin and Straus
\cite{MoSt65} following in turn from (\ref{maxmu1}) (see \cite{Nik06a}). The
implications
\[
\eqref{maxmu}\Longrightarrow\left(  \ref{maxmu1}\right)  \Longrightarrow
\text{MS}\Longrightarrow\eqref{maxmu}
\]
justify regarding inequality \eqref{maxmu} as a spectral form of Tur\'{a}n's
theorem, well suited for nontrivial generalizations. For example, the
following conjecture seems to be quite subtle.

\begin{conjecture}
Let $G$ be a $K_{r+1}$-free graph with $m$ edges. Then%
\[
\mu_{1}^{2}\left(  G\right)  +\mu_{2}^{2}\left(  G\right)  \leq\frac{r-1}%
{r}\text{ }2m.
\]

\end{conjecture}

If true, this conjecture is best possible whenever $r$ divides $n$. Indeed,
for $r|n$, $n=qr$, the Tur\'{a}n graph $T_{r}(n)$ (i.e., the complete
$r$-partite graph $K_{r}(q)$ with $q$ vertices in each class) has
$r(r-1)q^{2}/2$ edges, and there are three eigenvalues: $(r-1)q$, with
multiplicity $1$, $-q$, with multiplicity $r-1$, and $0$, with multiplicity
$r(q-1)$, so that $\mu_{1}(G)=(r-1)q$ and $\mu_{2}(G)=0$.\bigskip

The aim of this note is to prove further relations between $\mu\left(
G\right)  $ and the number of cliques in $G$. In \cite{Nik02} it is proved
that
\begin{equation}
\mu^{\omega}\left(  G\right)  \leq\sum_{s=2}^{\omega}\left(  s-1\right)
k_{s}\left(  G\right)  \mu^{\omega-s}\left(  G\right)  \label{polyn}%
\end{equation}
with equality holding if and only if $G$ is a complete $\omega$-partite graph
with possibly some isolated vertices. It turns out that this inequality is one
of a whole sequence of similar inequalities.

\begin{theorem}
\label{le3mu}For every graph $G$ and $r\geq2,$%
\[
\mu^{r+1}\left(  G\right)  \leq\left(  r+1\right)  k_{r+1}\left(  G\right)
+\sum_{s=2}^{r}\left(  s-1\right)  k_{s}\left(  G\right)  \mu^{r+1-s}\left(
G\right)  .
\]

\end{theorem}

Observe that, with $r=\omega+1$, Theorem \ref{le3mu} implies (\ref{polyn}).
Theorem \ref{le3mu} also implies a lower bound on the number of cliques of any
given order, as stated below.

\begin{theorem}
\label{tmomo}For every graph $G=G\left(  n\right)  $ and $r\geq2,$
\[
k_{r+1}\left(  G\right)  \geq\left(  \frac{\mu\left(  G\right)  }{n}%
-1+\frac{1}{r}\right)  \frac{r\left(  r-1\right)  }{r+1}\left(  \frac{n}%
{r}\right)  ^{r+1}.
\]

\end{theorem}

We also prove the following extension of an earlier result of
ours~\cite{BoNi04}.

\begin{theorem}
\label{leNSMM}Let $1\leq s\leq r<\omega\left(  G\right)  $ and $\alpha\geq0$.
If $G=G\left(  n\right)  $ and%
\begin{equation}
\left(  s+1\right)  k_{s+1}\left(  G\right)  \geq n^{s+1}\prod_{t=1}%
^{s}\left(  \frac{r-t}{rt}+\alpha\right)  , \label{cond1}%
\end{equation}
then%
\begin{equation}
k_{r+1}\left(  G\right)  \geq\alpha\frac{r^{2}}{r+1}\left(  \frac{n}%
{r}\right)  ^{r+1}. \label{lok}%
\end{equation}

\end{theorem}

Note that Theorems \ref{leNSMM} and \ref{tmomo} hold for all values of the
parameters satisfying the conditions there; in particular, $\alpha$ may depend
on $n$.

Our final theorem is the following stability result.

\begin{theorem}
\label{tstab}For all $r\geq2$ and $0\leq\alpha\leq2^{-10}r^{-6},$ if
$G=G\left(  n\right)  $ is a $K_{r+1}$-free graph with%
\begin{equation}
\mu\left(  G\right)  \geq\left(  1-\frac{1}{r}-\alpha\right)  n, \label{reqmu}%
\end{equation}
then $G$ contains an induced $r$-partite graph $G_{0}$ of order $v\left(
G_{0}\right)  >\left(  1-3\alpha^{1/3}\right)  n$ and minimum degree%
\[
\delta\left(  G_{0}\right)  >\left(  1-\frac{1}{r}-6\alpha^{1/3}\right)  n.
\]

\end{theorem}

\section{Proofs}

\subsection{Proof of Theorem \ref{le3mu}}

For a vertex $u\in V\left(  G\right)  $, write $w_{l}\left(  u\right)  $ for
the number of $l$-walks starting with $u$ and $k_{r}\left(  u\right)  $ for
the number of $r$-cliques containing $u.$ Clearly, it is enough to prove the
assertion for $2\leq r<\omega\left(  G\right)  $, since the case $r\geq
\omega\left(  G\right)  $ follows easily from (\ref{polyn}).

It is shown in \cite{Nik02} that for all $2\leq s\leq\omega\left(  G\right)  $
and $l\geq2,$%
\begin{equation}
\sum_{u\in V\left(  G\right)  }\big(k_{s}\left(  u\right)  w_{l+1}\left(
u\right)  -k_{s+1}\left(  u\right)  w_{l}\left(  u\right)  \big)\leq\left(
s-1\right)  k_{s}\left(  G\right)  w_{l}\left(  G\right)  . \label{oldin}%
\end{equation}
Summing these inequalities for $s=2,...r,$ we obtain%
\[
\sum_{u\in V\left(  G\right)  }\big(k_{2}\left(  u\right)  w_{l+r-1}\left(
u\right)  -k_{r+1}\left(  u\right)  w_{l}\left(  u\right)  \big)\leq\sum
_{s=2}^{r}\left(  s-1\right)  k_{s}\left(  G\right)  w_{l+r-s}\left(
G\right)  ,
\]
and so, after rearranging,%
\[
w_{l+r}\left(  G\right)  -\sum_{s=2}^{r}\left(  s-1\right)  k_{i}\left(
G\right)  w_{l+r-s}\left(  G\right)  \leq\sum_{u\in V\left(  G\right)
}k_{r+1}\left(  u\right)  w_{l}\left(  u\right)  .
\]
Noting that $w_{l}\left(  u\right)  \leq w_{l-1}\left(  G\right)  ,$ this
implies that
\[
\sum_{u\in V\left(  G\right)  }k_{r+1}\left(  u\right)  w_{l}\left(  u\right)
\leq w_{l-1}\left(  G\right)  \sum_{u\in V\left(  G\right)  }k_{r+1}\left(
u\right)  =\left(  r+1\right)  k_{r+1}\left(  G\right)  w_{l-1}\left(
G\right)  ,
\]
and so,
\[
\frac{w_{l+r}\left(  G\right)  }{w_{l-1}\left(  G\right)  }-\sum_{s=2}%
^{r}\left(  s-1\right)  k_{s}\left(  G\right)  \frac{w_{l+r-s}\left(
G\right)  }{w_{l-1}\left(  G\right)  }\leq\left(  r+1\right)  k_{r+1}\left(
G\right)  .
\]

Given $n$, there are non-negative constants $c_{1},\dots,c_{n}$ such that for
$G=G\left(  n\right)  $ we have
\[
w_{l}\left(  G\right)  =c_{1}\mu_{1}^{l-1}\left(  G\right)  +\cdots+c_{n}%
\mu_{n}^{l-1}\left(  G\right)  ,
\]
(See, e.g., \cite{CDS80}, p. 44.) Since $\omega>2$, our graph $G$ is not
bipartite and so $\left\vert \mu_{n}(G)\right\vert <\mu_{1}(G)$. Therefore,
for every fixed $q$, we have
\[
\lim_{l\rightarrow\infty}\frac{w_{l+q}\left(  G\right)  }{w_{l-1}\left(
G\right)  }=\mu^{q+1}\left(  G\right)  ,
\]
and the assertion follows. \hfill{$\Box$}

\subsection{Proof of Theorem \ref{leNSMM}}

Moon and Moser \cite{MoMo62} stated the following result (for a proof see
\cite{KhNi78} or \cite{Lov79}, Problem 11.8): if $G=G\left(  n\right)  $ and
$k_{s}\left(  G\right)  >0$, then
\[
\frac{\left(  s+1\right)  k_{s+1}\left(  G\right)  }{sk_{s}\left(  G\right)
}-\frac{n}{s}\geq\frac{sk_{s}\left(  G\right)  }{\left(  s-1\right)
k_{s-1}\left(  G\right)  }-\frac{n}{s-1}.
\]
Equivalently, for $1\leq s<t<\omega\left(  G\right)  $, we have%
\begin{equation}
\frac{\left(  t+1\right)  k_{t+1}\left(  G\right)  }{tk_{t}\left(  G\right)
}-\frac{n}{t}\geq\frac{\left(  s+1\right)  k_{s+1}\left(  G\right)  }%
{sk_{s}\left(  G\right)  }-\frac{n}{s}. \label{MoMo}%
\end{equation}
Let $s\in\left[  1,r\right]  $ be the smallest integer for which (\ref{cond1})
holds. This implies either $s=1$ or%
\begin{equation}
sk_{s}\left(  G\right)  <n^{s}\prod_{t=1}^{s-1}\left(  \frac{r-t}{rt}%
+\alpha\right)  \label{secin}%
\end{equation}
for some $s\in\left[  2,r\right]  $. Suppose first that $s=1$. (This case is
considered in \cite{BoNi04}, but for the sake of completeness we present it
here.) We have%
\[
\frac{2k_{2}\left(  G\right)  }{k_{1}\left(  G\right)  }-n\geq\left(
\frac{r-1}{r}+\alpha\right)  n-n=\alpha n-\frac{n}{r},
\]
and so, for all $t=1,\ldots,r$, inequality (\ref{MoMo}) implies that
\[
\frac{\left(  t+1\right)  k_{t+1}\left(  G\right)  }{tk_{t}\left(  G\right)
}\geq\alpha n+\frac{n}{t}-\frac{n}{r}.
\]
Multiplying these inequalities for $t=1,\ldots,r$, we obtain that%
\[
\left(  r+1\right)  k_{r+1}\left(  G\right)  \geq n^{r+1}\prod_{t=1}%
^{r}\left(  \frac{r-t}{rt}+\alpha\right)  \geq\alpha r^{2}\left(  \frac{n}%
{r}\right)  ^{r+1}\prod_{t=1}^{r-1}\frac{r-t}{t}=\alpha r^{2}\left(  \frac
{n}{r}\right)  ^{r+1},
\]
proving the result in this case.

Assume now that (\ref{secin}) holds for some $s\in\left[  2,r\right]  $. Then
we have%
\[
\frac{\left(  s+1\right)  k_{s+1}\left(  G\right)  }{sk_{s}\left(  G\right)
}>\left(  \frac{r-s}{rs}+\alpha\right)  n.
\]
and so, for every $t=s,...,r,$%
\[
\frac{\left(  t+1\right)  k_{t+1}\left(  G\right)  }{tk_{t}\left(  G\right)
}>\frac{n}{t}-\frac{n}{s}+\frac{r-s}{rs}n+\alpha n=\left(  \frac{r-t}%
{rt}+\alpha\right)  n.
\]
Multiplying these inequalities for $t=s+1,...,r,$ we obtain%
\[
\frac{\left(  r+1\right)  k_{r+1}\left(  G\right)  }{\left(  s+1\right)
k_{s+1}\left(  G\right)  }>n^{r-s}\prod_{t=s+1}^{r}\left(  \frac{r-t}%
{rt}+\alpha\right)  .
\]
Appealing to (\ref{cond1}), this implies that
\[
\left(  r+1\right)  k_{r+1}\left(  G\right)  >n^{r+1}\prod_{t=1}^{r}\left(
\frac{r-t}{rt}+\alpha\right)  =\alpha n^{r+1}\prod_{t=1}^{r-1}\left(
\frac{r-t}{rt}+\alpha\right)  \geq\alpha r^{2}\left(  \frac{n}{r}\right)
^{r+1},
\]
as required. \hfill{$\Box$}

\subsection{Proof of Theorem \ref{tmomo}}

Set
\[
\alpha=\frac{\mu}{n}-1+\frac{1}{r-1}.
\]
Clearly we may assume that $\alpha>0$, since otherwise the assertion is
trivial. Suppose that
\begin{equation}
sk_{s}\left(  G\right)  >n^{s}\prod_{t=1}^{s-1}\left(  \frac{r-t}{rt}%
+\alpha\right)  \label{in2}%
\end{equation}
for some $s\in\left[  2,r\right]  $. Then, by Theorem \ref{leNSMM},%
\[
\left(  r+1\right)  k_{r+1}\left(  G\right)  >\alpha\frac{r^{2}}{r+1}\left(
\frac{n}{r}\right)  ^{r+1}\geq\alpha\frac{r\left(  r-1\right)  }{r+1}\left(
\frac{n}{r}\right)  ^{r+1},
\]
completing the proof. Thus we may and shall assume that (\ref{in2}) fails for
every $s\in\left[  r-1\right]  $.

From Theorem \ref{le3mu} we have
\begin{equation}
\left(  r+1\right)  k_{r+1}\left(  G\right)  \geq\mu^{r+1}\left(  G\right)
-\sum_{s=2}^{r}\left(  s-1\right)  k_{s}\left(  G\right)  \mu^{r+1-s}\left(
G\right)  . \label{minK}%
\end{equation}
Substituting the bounds on $k_{s}\left(  G\right)  $ into (\ref{minK}), and
setting $\mu=\mu\left(  G\right)  /n$, we obtain
\begin{align*}
\frac{\left(  r+1\right)  k_{r+1}\left(  G\right)  }{n^{r+1}}  &  \geq
\mu^{r+1}-\sum_{s=2}^{r}\mu^{r+1-s}\frac{s-1}{s}\prod_{t=1}^{s-1}\left(
\frac{r-t}{rt}+\alpha\right) \\
&  \geq\mu^{r+1}-\mu^{r+1-2}\frac{1}{2}\left(  \frac{r-1}{r}+\alpha\right)
+\sum_{s=3}^{r}\frac{s-1}{s}\mu^{r+1-s}\prod_{t=1}^{s-1}\left(  \frac{r-t}%
{rt}+\alpha\right) \\
&  \geq\mu^{r+1-2}\left(  \mu^{2}-\frac{1}{2}\left(  \frac{r-1}{r}%
+\alpha\right)  \right)  +\sum_{s=3}^{r}\frac{s-1}{s}\mu^{r+1-s}\prod
_{t=1}^{s-1}\left(  \frac{r-t}{rt}+\alpha\right) \\
&  \geq\mu^{r+1-2}\left(  \frac{r-1}{r}+\alpha\right)  \left(  \frac{r-2}%
{2r}+\alpha\right)  +\sum_{s=3}^{r}\frac{s-1}{s}\mu^{r+1-s}\prod_{t=1}%
^{s-1}\left(  \frac{r-t}{rt}+\alpha\right)  .
\end{align*}
By induction on $k$ we prove that, for all $k=2,\ldots,r,$%
\[
\frac{\left(  r+1\right)  k_{r+1}\left(  G\right)  }{n^{r+1}}\geq\mu
^{r+1-k}\prod_{t=1}^{k}\left(  \frac{r-t}{rt}+\alpha\right)  -\sum_{s=k+1}%
^{r}\frac{s-1}{s}\mu^{r+1-s}\prod_{t=1}^{s-1}\left(  \frac{r-t}{rt}%
+\alpha\right)
\]
and hence,%
\[
\frac{\left(  r+1\right)  k_{r+1}\left(  G\right)  }{n^{r+1}}\geq\mu
\prod_{t=1}^{r}\left(  \frac{r-t}{rt}+\alpha\right)  \geq\alpha\frac{r-1}%
{r}\prod_{t=1}^{r-1}\frac{r-t}{rt}=\alpha\frac{r-1}{r^{r}}.
\]
It follows that%
\[
k_{r+1}\left(  G\right)  \geq\alpha\frac{r\left(  r-1\right)  }{r+1}\left(
\frac{n}{r}\right)  ^{r+1},
\]
as required. \hfill{$\Box$}

\subsection{Proof of Theorem \ref{tstab}}

Inequality (\ref{maxmu}) for $s=2$ together with (\ref{reqmu}) implies that%
\[
2\frac{r-1}{r}e\left(  G\right)  \geq\mu^{2}\left(  G\right)  >\left(
\frac{r-1}{r}-\alpha\right)  ^{2}n^{2}>\left(  \left(  \frac{r-1}{r}\right)
^{2}-2\alpha\frac{r-1}{r}\right)  n^{2},
\]
and so,%
\[
e\left(  G\right)  \geq\left(  \frac{r-1}{2r}-2\alpha\right)  n^{2}.
\]
To complete our proof, let us recall the following stability theorem proved by
Nikiforov and Rousseau in \cite{NiRo05}. Let $r\geq2$ and $0<\beta\leq
2^{-9}r^{-6}$, and let $G=G(n)$ be a $K_{r+1}$-free graph satisfying
\[
e\left(  G\right)  \geq\left(  \frac{r-1}{2r}-\beta\right)  n^{2}.
\]
Then $G$ contains an induced $r$-partite graph $G_{0}$ of order $v\left(
G_{0}\right)  >\left(  1-2\alpha^{1/3}\right)  n$ and with minimum degree
\[
\delta\left(  G_{0}\right)  \geq\left(  1-\frac{1}{r}-4\beta^{1/3}\right)  n.
\]
Setting $\beta=2\alpha,$ in view of $4\cdot2^{1/3}<6,$ the required
inequalities follow. \hfill{$\Box$}

\textbf{Acknowledgement}. Part of this research was completed while the
authors were enjoying the hospitality of the Institute for Mathematical
Sciences, National University of Singapore in 2006.

\end{document}